\documentclass[10pt]{article}
\begin{document}
\title{ {\bf A Characterization On Potentially $K_6-C_4$-graphic
Sequences}
\thanks{   Project Supported by  NSF of Fujian(2008J0209),
 Fujian Provincial Training
Foundation for "Bai-Quan-Wan Talents Engineering" , Project of
Fujian Education Department and Project of Zhangzhou Teachers
College.}}
\author{{ Lili Hu , Chunhui Lai}\\
{\small Department of Mathematics, Zhangzhou Teachers College,}
\\{\small Zhangzhou, Fujian 363000,
 P. R. of CHINA.}\\{\small  jackey2591924@163.com ( Lili Hu, Corresponding author)}
 \\{\small   zjlaichu@public.zzptt.fj.cn(Chunhui
 Lai )}
}

\date{}
\maketitle
\begin{center}
\begin{minipage}{4.1in}
\vskip 0.1in
\begin{center}{\bf Abstract}\end{center}
 { For given a graph $H$, a graphic sequence $\pi=(d_1,d_2,\cdots,d_n)$ is said to be potentially
 $H$-graphic if there exists a realization of $\pi$ containing $H$ as a subgraph.
 Let $K_m-H$ be the graph
obtained from $K_m$ by removing the edges set $E(H)$ where $H$ is a
subgraph of $K_m$. In this paper, we characterize the potentially
$K_6-C_4$-graphic sequences. This characterization implies a
 theorem due to  Hu and Lai [7].}
\par
\end{minipage}
\end{center}
 \par
 \section{Introduction}
\par
\baselineskip 14pt
    We consider finite simple graphs. Any undefined notation follows
that of Bondy and Murty $[1]$. The set of all non-increasing
nonnegative integer sequence $\pi=(d_1,d_2,\cdots,d_n)$ is denoted
by $NS_n$. A sequence $\pi\in NS_n$ is said to be graphic if it is
the degree sequence of a simple graph $G$ of order $n$; such a graph
$G$ is referred as a realization of $\pi$. The set of all graphic
sequences in $NS_n$ is denoted by $GS_n$.  Let $C_k$ and $P_k$
denote a cycle on $k$ vertices and a path on $k+1$ vertices,
respectively. Let $\sigma(\pi)$ be the sum of all the terms of
$\pi$, and let [x] be the largest integer less than or equal to $x$.
A graphic sequence $\pi$ is said to be potentially $H$-graphic if it
has a realization $G$ containing $H$ as a subgraph. Let $G-H$ denote
the graph obtained from $G$ by removing the edges set $E(H)$ where
$H$ is a subgraph of $G$.  In the degree sequence, $r^t$ means $r$
repeats $t$ times, that is, in the realization of the sequence there
are $t$ vertices of degree $r$.
\par

  Gould et al.[6] considered an extremal problem on potentially $H$-graphic sequences as follows:
  determine the smallest even
integer $\sigma(H,n)$ such that every n-term positive graphic
sequence $\pi$ with $\sigma(\pi)\geq \sigma(H,n)$ has a realization
$G$ containing $H$ as a subgraph. A harder question is to
characterize the potentially
 $H$-graphic sequences without zero terms. Yin and Li [18] gave two sufficient conditions
 for $\pi\in GS_n$ to be potentially $K_r-e$-graphic.  Luo [15] characterized the potentially
 $C_k$-graphic sequences for each $k=3,4,5$. Chen [2] characterized the potentially
 $C_6$-graphic sequences.  Chen et al.[3] characterized the potentially
 $_kC_l$-graphic sequences for each $3\leq k\leq5$, $l=6$. Recently, Luo and Warner [16] characterized the potentially
 $K_4$-graphic sequences.  Eschen and Niu [5] characterized the potentially
 $K_4-e$-graphic sequences. Yin et al.[19] characterized the potentially
 $_3C_4$, $_3C_5$ and $_4C_5$-graphic sequences. Yin and Chen [20] characterized the
 potentially $K_{r,s}$-graphic sequences for $r=2,s=3$ and
 $r=2,s=4$. Yin et al.[21] characterized the
 potentially $K_5-e$ and $K_6$-graphic sequences. Besides, Yin [22] characterized the potentially
$K_6-K_3$-graphic sequences. Chen and Li [4] characterized the
potentially $K_{1,t}+e$-graphic sequences. Xu and Lai[17]
characterized the potentially $K_6-C_5$-graphic sequences.  Hu and
Lai [8,11] characterized the potentially $K_5-C_4$ and
$K_5-E_3$-graphic sequences where $E_3$ denotes graphs with 5
vertices and 3 edges. In [12], they characterized the potentially
$K_{3,3}$ and $K_6-C_6$-graphic sequences. Moreover,
 Hu, Lai and Wang [10] characterized the potentially  $K_5-P_4$ and $K_5-Y_4$-graphic
sequences where $Y_4$ is a tree on 5 vertices and 3 leaves.

\par
In this paper, we characterize the potentially
$K_6-C_4$-graphic sequences. This characterization implies a
 theorem due to  Hu and Lai [7]. Up to now, the problem
of characterizing the potentially $K_6-C_k(3\leq k\leq6)$-graphic
sequences has been completely solved.
\par
\section{Preparations}\par
   Let $\pi=(d_1,\cdots,d_n)\in NS_n,1\leq k\leq n$. Let
    $$ \pi_k^{\prime\prime}=\left\{
    \begin{array}{ll}(d_1-1,\cdots,d_{k-1}-1,d_{k+1}-1,
    \cdots,d_{d_k+1}-1,d_{d_k+2},\cdots,d_n), \\ \mbox{ if $d_k\geq k,$}\\
    (d_1-1,\cdots,d_{d_k}-1,d_{d_k+1},\cdots,d_{k-1},d_{k+1},\cdots,d_n),
     \\ \mbox{if $d_k < k.$} \end{array} \right. $$
  Denote
  $\pi_k^\prime=(d_1^\prime,d_2^\prime,\cdots,d_{n-1}^\prime)$, where
  $d_1^\prime\geq d_2^\prime\geq\cdots\geq d_{n-1}^\prime$ is a
  rearrangement of the $n-1$ terms of $\pi_k^{\prime\prime}$. Then
  $\pi_k^{\prime}$ is called the residual sequence obtained by
  laying off $d_k$ from $\pi$. For simplicity, we denote $\pi_n^\prime$ by $\pi^\prime$ in this paper.
  We need the following results.
\par
    {\bf Theorem 2.1 [6]} If $\pi=(d_1,d_2,\cdots,d_n)$ is a graphic
 sequence with a realization $G$ containing $H$ as a subgraph,
 then there exists a realization $G^\prime$ of $\pi$ containing $H$ as a
 subgraph so that the vertices of $H$ have the largest degrees of
 $\pi$.\par
 \par
    {\bf Theorem 2.2 [8]} Let $\pi=(d_1,d_2,\cdots,d_n)\in NS_n$ be a graphic sequence with $n\geq5$.
Then $\pi$ is potentially $K_5-C_4$-graphic if and only if the
following conditions hold: \par (1) $d_1\geq4$, $d_5\geq2$;
\par
(2)$\pi\neq(4,2^5), (4,2^6), ((n-2)^2,2^{n-2}),
(n-k,k+i,2^i,1^{n-i-2})$ where $i=3,4,\cdots,n-2k$ and
$k=1,2,\cdots,[\frac{n-1}{2}]-1$.

\par
    {\bf Lemma 2.3 [9]} If $\pi=(d_1,d_2,\cdots,d_n)$ is a
 nonincreasing sequence of positive integers with even $\sigma(\pi)$,  $n\geq4$,
$d_1\leq3$ and $\pi\neq(3^3,1),(3^2,1^2)$, then $\pi$ is graphic.

\par
    {\bf Lemma 2.4 [22]}  Let $\pi=(4^x,3^y,2^z,1^m)$ with even $\sigma(\pi)$, $x+y+z+m=n\geq5$
    and $x\geq1$.
Then $\pi\in GS_n$ if and only if $\pi \not \in A$, where
$A=\{(4,3^2,1^2), (4,3,1^3), (4^2,2,1^2), (4^2,3,2,1), (4^3,1^2),
(4^3,2^2), (4^3,3,1), (4^4,2)$,\par
 $(4^2,3,1^3), (4^2,1^4), (4^3,2,1^2), (4^4,1^2), (4^3,1^4)\}$.\par

    {\bf Lemma 2.5 (Kleitman and Wang [13])}\ \   $\pi$ is
graphic if and only if $\pi_k^\prime$ is graphic.
 \par
    The following corollary is obvious.\par
\par
    {\bf Corollary 2.6}\ \    Let $H$ be a simple graph. If $\pi^\prime$ is
 potentially $H$-graphic, then $\pi$ is
 potentially $H$-graphic.

\par
\section{ Main Theorems} \par
\par
\textbf{\noindent Theorem 3.1}  Let $\pi=(d_1,d_2,\cdots,d_n)\in
NS_n$ be a graphic sequence with $n\geq6$. Then $\pi$ is potentially
$K_6-C_4$-graphic if and only if the following conditions hold:
\par
  $(1)$ $d_2\geq5$, $d_6\geq3$;\par
  $(2)$ $\pi=(d_1,d_2,d_3,3^k,2^t,1^{n-3-k-t})$  implies $d_1+d_2+d_3\leq n+2k+t+1$;\par
  $(3)$ $\pi\neq(5^2,4^6)$,  $(5^2,4^7)$,  $(6^2,3^6)$,
  $(6,5,4,3^5)$, $(6,5,3^7)$,  $(5^3,4,3^3)$,  $(5^3,3^5)$,  $(5^2,4^2,3^4)$,
  $(5^2,4,3^6)$, $(5^2,4,3^4)$,  $(5^2,3^6)$,  $(6,5,3^5,2)$,
  $(5^3,3^3,2)$, $(5^2,4,3^4,2)$, $(5^2,3^6,2)$, $(5^2,3^4,2)$,
  $(5^2,3^4,2^2)$, $(6,5,3^6,1)$, $(5^3,3^4,1)$, $(5^2,4,3^5,1)$,
  $(5^2,3^7,1)$,  $(5^2,3^6,1^2)$,  $(5^2,3^5,1)$,
   $(n-1,5,3^5,1^{n-7})$, $(n-1,5,3^6,1^{n-8})$.
\par
{\bf Proof:} First we show the conditions (1)-(3) are necessary
conditions for $\pi$ to be potentially $K_6-C_4$-graphic. Assume
that $\pi$ is potentially $K_6-C_4$-graphic. $(1)$ is obvious.  If
$\pi=(d_1,d_2,d_3,3^k,2^t,1^{n-3-k-t})$ is potentially
$K_6-C_4$-graphic, then according to Theorem 2.1, there exists a
realization $G$ of $\pi$ containing $K_6-C_4$ as a subgraph so that
the vertices of $K_6-C_4$ have the largest degrees of $\pi$.
Therefore, the sequence
$\pi_1=(d_1-5,d_2-5,d_3-3,3^{k-3},2^t,1^{n-3-k-t})$ obtained from
$G-(K_6-C_4)$ is graphic and there exists no edge among three
vertices with degree $d_1-5$, $d_2-5$ and $d_3-3$ in the realization
of $\pi_1$. It follows $d_1-5+d_2-5+d_3-3\leq3(k-3)+2t+n-3-k-t$,
i.e., $d_1+d_2+d_3\leq n+2k+t+1$.  Hence, (2) holds. Now it is easy
to check that $(5^2,4^6)$, $(5^2,4^7)$, $(6^2,3^6)$, $(6,5,4,3^5)$,
$(6,5,3^7)$, $(5^3,4,3^3)$, $(5^3,3^5)$, $(5^2,4^2,3^4)$,
$(5^2,4,3^6)$, $(5^2,4,3^4)$, $(5^2,3^6)$, $(6,5,3^5,2)$,
$(5^3,3^3,2)$, $(5^2,4,3^4,2)$, $(5^2,3^6,2)$, $(5^2,3^4,2)$,
$(5^2,3^4,2^2)$, $(6,5,3^6,1)$, $(5^3,3^4,1)$, $(5^2,4,3^5,1)$,
$(5^2,3^7,1)$, $(5^2,3^6,1^2)$ and $(5^2,3^5,1)$ are not potentially
$K_6-C_4$-graphic. Since $(4,2^5)$ and $(4,2^6)$ are not potentially
$K_5-C_4$-graphic by Theorem 2.2, we have
$\pi\neq(n-1,5,3^5,1^{n-7})$ and $(n-1,5,3^6,1^{n-8})$.  Hence, (3)
holds.
\par
 To prove the sufficiency, we use induction on $n$. Suppose the graphic sequence
$\pi$ satisfies the conditions (1)-(3). We first prove the base case
where $n=6$. Since $\pi\neq(5^3,3^3)$, then $\pi$ is one of the
following: $(5^6)$, $(5^4,4^2)$, $(5^3,4^2,3)$, $(5^2,4^4)$,
 $(5^2,4^2,3^2)$, $(5^2,3^4)$. It is
easy to check that all of these are potentially $K_6-C_4$-graphic.
Now suppose that the sufficiency holds for $n-1(n\geq7)$, we will
show that $\pi$ is potentially $K_6-C_4$-graphic in terms of the
following cases:
\par
\textbf{Case 1:} $d_n\geq4$.  Consider
$\pi^\prime=(d_1^\prime,d_2^\prime,\cdots,d_{n-1}^\prime)$ where
$d_{n-3}^\prime\geq4$ and $d_{n-1}^\prime\geq3$. If $\pi^\prime$
satisfies $(1)$ and $(3)$, then by the induction hypothesis,
$\pi^\prime$ is potentially $K_6-C_4$-graphic, and hence so is
$\pi$.
\par
  If $\pi^\prime$ does not satisfy $(1)$, i.e., $d_2^\prime=4$, then
  $d_2=5$. We will proceed with the following two cases $d_1=5$ and
  $d_1\geq6$.
\par
\textbf{Subcase 1:} $d_1=5$. Then $\pi=(5^k,4^{n-k})$ where $1\leq
k\leq5$. Since $\sigma(\pi)$ is even, we have $k=2$ or $k=4$. If
$k=2$, then $\pi=(5^2,4^{n-2})$. Since $\pi\neq(5^2,4^6)$ and
$(5^2,4^7)$, we have $n=7$ or $n\geq10$. It is easy to check  that
$(5^2,4^5)$ and $(5^2,4^8)$ are potentially $K_6-C_4$-graphic. If
$n\geq11$, let $\pi_1=(5^2,4^4)$, $\pi_2=(4^{n-6})$. Then by lemma
2.4, $\pi_2$ is graphic. Let $G_1$ be a realization of $\pi_2$, then
$(K_6-2K_2)\cup G_1$ is a realization of $\pi$. Thus,
$\pi=(5^2,4^{n-2})(n\neq8,9)$ is potentially $K_6-C_4$-graphic since
$K_6-C_4\subseteq K_6-2K_2$. Similarly, one can show that
$\pi=(5^4,4^{n-4})$ is also potentially $K_6-C_4$-graphic.
\par
\textbf{Subcase 2:} $d_1\geq6$. Then $\pi=(d_1,5^k,4^{n-1-k})$ where
$1\leq k\leq3$, $d_1$ and $k$ have the same parity. We will show
that $\pi$ is potentially $K_6-C_4$-graphic.
\par
If $k=1$, then $\pi=(d_1,5,4^{n-2})$ where $d_1$ is odd. If
$n\leq10$, then $\pi$ is one of the following: $(7,5,4^6)$,
$(7,5,4^7)$, $(7,5,4^8)$, $(9,5,4^8)$. It is easy to check that all
of these are potentially $K_6-C_4$-graphic. If $n\geq11$, let
$\pi_1=(5^2,4^4)$, $\pi_2=(d_1-5,4^{n-6})$. Then the residual
sequence $\pi_2^\prime=(4^{n-1-d_1},3^{d_1-5})$ obtained by laying
off $d_1-5$ from $\pi_2$ is graphic by lemma 2.3 and lemma 2.4.
Hence, $\pi_2$ is graphic. Let $G_1$ be a realization of $\pi_2$,
and $x\in V(G_1)$ with $d_{G_1}(x)=d_1-5$. Denote $G=(K_{1,2,2}\cup
G_1)\cup \{xx_1, xx_2, xx_3, xx_4, xx_5\}$ where $x_i\in
V(K_{1,2,2}), i=1,\cdots,5.$ i.e., $G$ is the graph obtained from
$K_{1,2,2}\cup G_1$ by adding new edges $xx_1, xx_2, xx_3, xx_4,
xx_5$ to $K_{1,2,2}\cup G_1$. Clearly, $G$ is a realization of $\pi$
and contains $K_6-C_4$.
\par
If $k=2$, then $\pi=(d_1,5^2,4^{n-3})$ where $d_1$ is even. If
$n\leq11$, then $\pi$ is one of the following: $(6,5^2,4^4)$,
$(6,5^2,4^5)$, $(6,5^2,4^6)$, $(6,5^2,4^7)$, $(6,5^2,4^8)$,
$(8,5^2,4^6)$, $(8,5^2,4^7)$, $(8,5^2,4^8)$, $(10,5^2,4^8)$. It is
easy to check that all of these are potentially $K_6-C_4$-graphic.
If $n\geq12$, let $\pi_1=(6,5^2,4^4)$, $\pi_2=(d_1-6,4^{n-7})$. Then
the residual sequence $\pi_2^\prime=(4^{n-1-d_1},3^{d_1-6})$
obtained by laying off $d_1-6$ from $\pi_2$ is graphic by lemma 2.3
and lemma 2.4. Hence, $\pi_2$ is graphic. Let $G_1$ be a realization
of $\pi_2$ and $x\in V(G_1)$ with $d_{G_1}(x)=d_1-6$. Denote
$G=(K_6-P_5\cup G_1)\cup \{xx_1, xx_2, xx_3, xx_4, xx_5, xx_6\}$
where $x_i\in V(K_6-P_5), i=1,\cdots,6.$ i.e., $G$ is the graph
obtained from $K_6-P_5\cup G_1$ by adding new edges $xx_1, xx_2,
xx_3, xx_4, xx_5, xx_6$ to $K_6-P_5\cup G_1$. Clearly, $G$ is a
realization of $\pi$ and contains $K_6-C_4$.
\par
Similarly, one can
show that $\pi=(d_1,5^3,4^{n-4})$ is also potentially
$K_6-C_4$-graphic.
\par
If $\pi^\prime$ does not satisfy $(3)$, then $\pi^\prime=(5^2,4^6)$
or $(5^2,4^7)$. Hence, $\pi=(6^2,5^2,4^5)$, $(6,5^4,4^4)$,
$(5^6,4^3)$, $(6^2,5^2,4^6)$, $(6,5^4,4^5)$, $(5^6,4^4)$. It is easy
to check that all of these are potentially $K_6-C_4$-graphic.
\par
\textbf{Case 2:} $d_n=3$. Consider
$\pi^\prime=(d_1^\prime,d_2^\prime,\cdots,d_{n-1}^\prime)$ where
$d_2^\prime\geq4$ and $d_{n-2}^\prime\geq3$. If $\pi^\prime$
satisfies $(1)$-$(3)$, then by the induction hypothesis,
$\pi^\prime$ is potentially $K_6-C_4$-graphic, and hence so is
$\pi$.
\par
  If $\pi^\prime$ does not satisfy $(1)$, there are three subcases:
  \par
\textbf{Subcase 1:} $d_2^\prime\geq5$ and $d_6^\prime=2$. Then
$\pi=(d_1,d_2,3^5)$ where $d_1\geq d_2\geq6$, which is impossible.
\par
\textbf{Subcase 2:} $d_2^\prime=4$ and $d_6^\prime=2$. Then
$\pi=(6,5,3^5)$, which contradicts condition (3).
\par
\textbf{Subcase 3:} $d_2^\prime=4$ and $d_6^\prime\geq3$. Then
$d_2=5$, we will proceed with the following two cases $d_1=5$ and
$d_1\geq6$.
\par
\textbf{Subcase 3.1:} $d_1=5$. Then $\pi=(5^k,4^t,3^{n-k-t})$ where
$2\leq k\leq4$, $n-k-t\geq1$ and $n-t$ is even.
\par
If $k=2$, then $\pi=(5^2,4^t,3^{n-2-t})$. We will show that $\pi$ is
potentially $K_6-C_4$-graphic. If $t=0$, then $\pi=(5^2,3^{n-2})$.
Since $\pi\neq(5^2,3^6)$, we have $n\geq10$. It is enough to show
$\pi_1=(3^{n-6})$ is graphic. It follows by lemma 2.3. If $t=1$,
then $\pi=(5^2,4,3^{n-3})$. Since $\pi\neq(5^2,4,3^4)$ and
$(5^2,4,3^6)$, we have $n\geq11$. We only need to show that
$\pi_1=(3^{n-6},1)$ is graphic. It also follows by lemma 2.3. If
$t=2$, then $\pi=(5^2,4^2,3^{n-4})$. Since $\pi\neq(5^2,4^2,3^4)$,
we have $n\geq10$. It is enough to show $\pi_1=(3^{n-6},1^2)$ is
graphic. It follows by lemma 2.3. If $t=3$, then
$\pi=(5^2,4^3,3^{n-5})$. Since $\pi_1=(3^{n-6},1^3)$ is graphic by
lemma 2.3, $\pi=(5^2,4^3,3^{n-5})$ is potentially $K_6-C_4$-graphic.
If $t\geq4$, let $\pi_1=(5^2,4^4)$, $\pi_2=(4^{t-4},3^{n-2-t})$. If
$n\geq11$, then by lemma 2.3 and lemma 2.4, $\pi_2$ is graphic. Let
$G_1$ be a realization of $\pi_2$, then $(K_6-2K_2)\cup G_1$ is a
realization of $\pi$. Since $K_6-C_4\subseteq K_6-2K_2$, $\pi$ is
potentially $K_6-C_4$-graphic. If $n\leq10$, then
$\pi=(5^2,4^4,3^2)$, $(5^2,4^4,3^4)$, $(5^2,4^5,3^2)$ or
$(5^2,4^6,3^2)$. It is easy to check that all of these are
potentially $K_6-C_4$-graphic.
\par
If $k=3$, then $\pi=(5^3,4^t,3^{n-3-t})$. We will show that $\pi$ is
potentially $K_6-C_4$-graphic. If $t=0$, then $\pi=(5^3,3^{n-3})$.
Since $\pi\neq(5^3,3^5)$, we have $n\geq10$. It is easy to check
that $(5^3,3^7)$ and $(5^3,3^9)$ are potentially $K_6-C_4$-graphic.
If $n\geq14$, let $\pi_1=(5^3,3^7)$, $\pi_2=(3^{n-10})$ and $G_1$ be
a realization of $\pi_1$ which contains  $K_6-C_4$. Then by lemma
2.3, $\pi_2$ is graphic. Let $G_2$ be a realization of $\pi_2$, then
$G_1\cup G_2$ is a realization of $\pi=(5^3,3^{n-3})$. Similarly,
one can show that $\pi=(5^3,4^t,3^{n-3-t})$ is potentially
$K_6-C_4$-graphic for the cases $t=1$ and $t=2$. If $t\geq3$, let
$\pi_1=(5^3,4^2,3)$, $\pi_2=(4^{t-2},3^{n-4-t})$. If $n\geq11$, then
$\pi_2$ is graphic by lemma 2.4. Let $G_1$ be a realization of
$\pi_2$, then $(K_6-P_2)\cup G_1$ is a realization of $\pi$. Since
$K_6-C_4\subseteq K_6-P_2$, $\pi$ is potentially $K_6-C_4$-graphic.
If $n\leq10$, then $\pi=(5^3,4^3,3)$, $(5^3,4^3,3^3)$,
$(5^3,4^4,3)$, $(5^3,4^4,3^3)$, $(5^3,4^5,3)$ or $(5^3,4^6,3)$. It
is easy to check that all of these are potentially
$K_6-C_4$-graphic.
\par
If $k=4$, then $\pi=(5^4,4^t,3^{n-4-t})$. We will show that $\pi$ is
potentially $K_6-C_4$-graphic. If $t=0$, then $\pi=(5^4,3^{n-4})$.
It is easy to check that $(5^4,3^4)$ and $(5^4,3^6)$ are potentially
$K_6-C_4$-graphic. If $n\geq12$, let $\pi_1=(5^4,3^4)$,
$\pi_2=(3^{n-8})$ and $G_1$ be a realization of $\pi_1$ which
contains  $K_6-C_4$. Then by lemma 2.3, $\pi_2$ is graphic. Let
$G_2$ be a realization of $\pi_2$, then $G_1\cup G_2$ is a
realization of $\pi=(5^4,3^{n-4})$. Similarly, one can show that
$\pi=(5^4,4^t,3^{n-4-t})$ is potentially $K_6-C_4$-graphic for the
cases $t=1$ and $t=2$. If $t\geq3$, let $\pi_1=(5^4,4^2)$,
$\pi_2=(4^{t-2},3^{n-4-t})$. If $n\geq11$, then $\pi_2$ is graphic
by lemma 2.4. Let $G_1$ be a realization of $\pi_2$, then
$(K_6-e)\cup G_1$ is a realization of $\pi$. Since $K_6-C_4\subseteq
K_6-e$, $\pi$ is potentially $K_6-C_4$-graphic. If $n\leq10$, then
$\pi=(5^4,4^3,3^2)$ or $(5^4,4^4,3^2)$. It is easy to check that
both of them are potentially $K_6-C_4$-graphic.
\par
\textbf{Subcase 3.2:} $d_1\geq6$. Then $\pi=(d_1,5,4^k,3^{n-2-k})$
where $n-2-k\geq1$, $d_1$ and $n-1-k$ have the same parity. We will
show that $\pi$ is potentially $K_6-C_4$-graphic.
\par
If $k=0$, then $\pi=(d_1,5,3^{n-2})$. Since $\pi\neq(6,5,3^5)$ and
$(7,5,3^6)$, we have $n\geq9$. If $n=9$, since $\pi\neq(6,5,3^7)$,
then $\pi=(8,5,3^7)$ which is potentially $K_6-C_4$-graphic. If
$n\geq10$, we only need to show that $\pi_1=(d_1-5,3^{n-6})$ is
graphic. Since the residual sequence
$\pi_1^\prime=(3^{n-1-d_1},2^{d_1-5})$ obtained by laying off
$d_1-5$ from $\pi_1$ is graphic by lemma 2.3, $\pi_1$ is graphic.
\par
If $k=1$, then $\pi=(d_1,5,4,3^{n-3})$. Since $\pi\neq(6,5,4,3^5)$ ,
we have $n\geq9$. It is enough to show $\pi_1=(d_1-5,3^{n-6},1)$ is
graphic and there exists no edge between two vertices with degree
$d_1-5$ and 1 in the realization of $\pi_1$. Hence, it suffices to
show $\pi_2=(3^{n-1-d_1},2^{d_1-5},1)$ is graphic. It follows by
lemma 2.3. With the same argument as above, one can show that
$\pi=(d_1,5,4^k,3^{n-2-k})$ is potentially $K_6-C_4$-graphic for the
cases $k=2$ and $k=3$.
\par
Now we consider the case where $k\geq4$. If $n\leq10$, then $\pi$ is
one of the following: $(6,5,4^4,3)$, $(6,5,4^4,3^3)$, $(6,5,4^5,3)$,
$(6,5,4^5,3^3)$, $(6,5,4^6,3)$, $(6,5,4^7,3)$, $(7,5,4^4,3^2)$,
$(7,5,4^4,3^4)$, $(7,5,4^5,3^2)$, $(7,5,4^6,3^2)$, $(8,5,4^4,3^3)$,
$(8,5,4^5,3^3)$, $(8,5,4^6,3)$, $(8,5,4^7,3)$, $(9,5,4^4,3^4)$,
$(9,5,4^6,3^2)$. It is easy to check that all of these are
potentially $K_6-C_4$-graphic. If $n\geq11$, let $\pi_1=(5^2,4^4)$,
$\pi_2=(d_1-5,4^{k-4},3^{n-2-k})$. Then the residual sequence
$\pi_2^\prime$ obtained by laying off $d_1-5$ from $\pi_2$ is
graphic by lemma 2.3 and lemma 2.4, and hence $\pi_2$ is also
graphic. Let $G_1$ be a realization of $\pi_2$ and $x\in V(G_1)$
with $d_{G_1}(x)=d_1-5$. Denote $G=(K_{1,2,2}\cup G_1)\cup \{xx_1,
xx_2, xx_3, xx_4, xx_5\}$ where $x_i\in V(K_{1,2,2}), i=1,\cdots,5.$
i.e., $G$ is the graph obtained from $K_{1,2,2}\cup G_1$ by adding
new edges $xx_1, xx_2, xx_3, xx_4, xx_5$ to $K_{1,2,2}\cup G_1$.
Clearly, $G$ is a realization of $\pi$ and contains $K_6-C_4$.
\par
If $\pi^\prime$ does not satisfy (2), there are two subcases:
\par
\textbf{Subcase 1:}
$\pi^\prime=(d_1^\prime,d_2^\prime,d_3^\prime,3^{n-4})$ and
$d_1^\prime+d_2^\prime+d_3^\prime> n-1+2(n-4)+1$, i.e.,
$d_1^\prime+d_2^\prime+d_3^\prime>3n-8$. If $d_3^\prime\leq 4$, then
$d_1^\prime+d_2^\prime> 3n-12$. It follows $3n-10\leq
d_1^\prime+d_2^\prime\leq2(n-2)$, i.e., $n\leq6$, a contradiction.
Thus, $d_3^\prime\geq 5$. Therefore, $\pi=(d_1,d_2,d_3,3^{n-3})$ and
$d_1+d_2+d_3> 3n-5$, a contradiction.
\par
\textbf{Subcase 2:} $\pi^\prime=(d_1^\prime,d_2^\prime,3^{n-4},2)$
and $d_1^\prime+d_2^\prime+3> n-1+2(n-5)+1+1$, i.e.,
$d_1^\prime+d_2^\prime>3n-12$. Hence, $3n-10\leq
d_1^\prime+d_2^\prime\leq2(n-2)$, i.e., $n\leq6$, a contradiction.
\par
If $\pi^\prime$ does not satisfy (3), since $\pi\neq(6^2,3^6)$, then
$\pi^\prime$ is one of the following: $(5^2,4^6)$, $(5^2,4^7)$,
$(6^2,3^6)$, $(6,5,4,3^5)$, $(6,5,3^7)$, $(5^3,4,3^3)$, $(5^3,3^5)$,
$(5^2,4^2,3^4)$, $(5^2,4,3^6)$, $(5^2,4,3^4)$, $(5^2,3^6)$,
$(6,5,3^5,2)$, $(5^2,3^6,2)$, $(6,5,3^5)$, $(7,5,3^6)$. Hence, $\pi$
is one of the following: $(6^2,5,4^5,3)$, $(6,5^3,4^4,3)$,
$(5^5,4^3,3)$, $(6^2,5,4^6,3)$,\ \  $(6,5^3,4^5,3)$,\ \
$(5^5,4^4,3)$,\ \ $(7^2,4,3^6)$,\ \  $(7,6,5,3^6)$,\ \
$(7,6,4^2,3^5)$,\ \ $(7,6,4,3^7)$, $(6^3,4,3^4)$, $(6^2,5^2,3^4)$,
$(6^3,3^6)$, $(6^2,5,4,3^5)$, $(6^2,4^3,3^4)$, $(6,5^3,3^5)$,
$(6^2,5,3^7)$, $(6^2,4^2,3^6)$, $(6^2,5,3^5)$, $(6^2,4^2,3^4)$,
$(6^2,4,3^6)$, $(7,6,3^7)$, $(6^2,3^8)$, $(7,6,4,3^5)$,
$(8,6,4,3^6)$. It is easy to check that all of these are potentially
$K_6-C_4$-graphic.
\par
\textbf{Case 3:} $d_n=2$. Consider
$\pi^\prime=(d_1^\prime,d_2^\prime,\cdots,d_{n-1}^\prime)$ where
$d_2^\prime\geq4$, $d_6^\prime\geq3$ and $d_{n-1}^\prime\geq2$. If
$\pi^\prime$
 satisfies $(1)$-$(3)$, then by the induction hypothesis,
$\pi^\prime$ is potentially $K_6-C_4$-graphic, and hence so is
$\pi$.
\par
  If $\pi^\prime$ does not satisfy $(1)$, i.e., $d_2^\prime=4$, then
$d_2=5$. There are two subcases:
\par
\textbf{Subcase 1:} $d_1\geq6$. Then
$\pi=(d_1,5,4^k,3^t,2^{n-2-k-t})$ where $k+t\geq4$, $n-2-k-t\geq1$,
and, $d_1$ and $t$ have different parities. We will show that $\pi$
is potentially $K_6-C_4$-graphic.
\par
If $k=0$, then $\pi=(d_1,5,3^t,2^{n-2-t})$. If $n\geq10$, we only
need to show that $\pi_1=(d_1-5,3^{t-4},2^{n-2-t})$ is graphic. The
residual sequence $\pi_1^\prime$ obtained by laying off $d_1-5$ from
$\pi_1$ clearly satisfies the hypothesis of lemma 2.3, and so
$\pi_1^\prime$ is graphic and hence so is $\pi_1$. If $n\leq9$,
since $\pi\neq(6,5,3^5,2)$, then $\pi=(6,5,3^5,2^2)$,
$(7,5,3^4,2^2)$, $(7,5,3^4,2^3)$, $(7,5,3^6,2)$, $(8,5,3^5,2^2)$. It
is easy to check that all of these are potentially
$K_6-C_4$-graphic.
\par
If $k=1$, then $\pi=(d_1,5,4,3^t,2^{n-3-t})$. If $n\geq9$, we only
need to show that $\pi_1=(d_1-5,3^{t-3},2^{n-3-t},1)$ is graphic and
there exists no edge between two vertices with degree $d_1-5$ and 1
in the realization of $\pi_1$. The residual sequence $\pi_1^\prime$
obtained by laying off $d_1-5$ from $\pi_1$ clearly satisfies the
hypothesis of lemma 2.3, and so $\pi_1^\prime$ is graphic and hence
so is $\pi_1$. If $n\leq8$, then $\pi=(6,5,4,3^3,2)$,
$(6,5,4,3^3,2^2)$, $(7,5,4,3^4,2)$. It is easy to check that all of
these are potentially $K_6-C_4$-graphic. With the same argument as
above, one can show that $\pi=(d_1,5,4^k,3^t,2^{n-2-k-t})$ is
potentially $K_6-C_4$-graphic for the cases $k=2$ and $k=3$.
\par
Now we consider the case where $k\geq4$. If $n\geq11$, let
$\pi_1=(5^2,4^4)$, $\pi_2=(d_1-5,4^{k-4},3^t,2^{n-2-k-t})$. Then the
residual sequence $\pi_2^\prime$ obtained by laying off $d_1-5$ from
$\pi_2$ is graphic by lemma 2.3 and lemma 2.4, and hence $\pi_2$ is
also graphic. Let $G_1$ be a realization of $\pi_2$ and and $x\in
V(G_1)$ with $d_{G_1}(x)=d_1-5$. Denote $G=(K_{1,2,2}\cup G_1)\cup
\{xx_1, xx_2, xx_3, xx_4, xx_5\}$ where $x_i\in V(K_{1,2,2}),
i=1,\cdots,5.$ i.e., $G$ is the graph obtained from $K_{1,2,2}\cup
G_1$ by adding new edges $xx_1, xx_2, xx_3, xx_4, xx_5$ to
$K_{1,2,2}\cup G_1$. Clearly, $G$ is a realization of $\pi$ and
contains $K_6-C_4$. If $n\leq10$, then $\pi$ is one of the
following: $(6,5,4^4,3,2)$, $(6,5,4^4,3,2^2)$, $(6,5,4^4,3,2^3)$,
$(6,5,4^4,3^3,2)$, $(6,5,4^5,3,2)$,\ \  $(6,5,4^5,3,2^2)$,\ \
$(6,5,4^6,3,2)$,\ \  $(7,5,4^4,2^2)$,\ \  $(7,5,4^4,2^3)$,\ \
$(7,5,4^4,2^4)$, \  $(7,5,4^4,3^2,2)$, \   $(7,5,4^4,3^2,2^2)$, \
$(7,5,4^5,2)$, \  $(7,5,4^5,2^2)$, \  $(7,5,4^5,2^3)$, \ \
$(7,5,4^5,3^2,2)$, \ \  $(7,5,4^6,2)$, \ \  $(7,5,4^6,2^2)$, \ \
$(7,5,4^7,2)$, \  \ $(8,5,4^4,3,2^2)$, $(8,5,4^4,3,2^3)$,
$(8,5,4^4,3^3,2)$, $(8,5,4^5,3,2)$, $(8,5,4^5,3,2^2)$,
$(8,5,4^6,3,2)$, $(9,5,4^4,2^4)$, \ $(9,5,4^4,3^2,2^2)$, \
$(9,5,4^5,2^3)$, \ $(9,5,4^5,3^2,2)$, \ $(9,5,4^6,2^2)$, \
$(9,5,4^7,2)$. It is easy to check that all of these are potentially
$K_6-C_4$-graphic.
\par
\textbf{Subcase 2:} $d_1=5$. Then $\pi=(5^i,4^k,3^t,2^{n-i-k-t})$
where $2\leq i\leq3$, $i+k+t\geq6$, $n-i-k-t\geq1$ and $i+t$ is
even. We will show that $\pi$ is potentially $K_6-C_4$-graphic.
\par
\textbf{Subcase 2.1:} $i=2$. Then $\pi=(5^2,4^k,3^t,2^{n-2-k-t})$.
If $k=0$, then $\pi=(5^2,3^t,2^{n-2-t})$. If $n\geq10$, it is enough
to show $\pi_1=(3^{t-4},2^{n-2-t})$ is graphic. It follows by lemma
2.3. If $n\leq9$, since $\pi\neq(5^2,3^4,2)$, $(5^2,3^4,2^2)$ and
$(5^2,3^6,2)$, we have $\pi=(5^2,3^4,2^3)$ which is potentially
$K_6-C_4$-graphic.
\par
If $k=1$, then $\pi=(5^2,4,3^t,2^{n-3-t})$. Since
$\pi\neq(5^2,4,3^4,2)$, we have $n\geq9$. It is enough to show
$\pi_1=(3^{t-3},2^{n-3-t},1)$ is graphic. It follows by lemma 2.3.
With the same argument as above, one can show that
$\pi=(5^2,4^k,3^t,2^{n-2-k-t})$ is potentially $K_6-C_4$-graphic for
the cases $k=2$ and $k=3$.
\par
Now we consider the case where $k\geq4$. If $n\geq11$, let
$\pi_1=(5^2,4^4)$, $\pi_2=(4^{k-4},3^t,2^{n-2-k-t})$. If
$\pi_2\neq(4^3,2^2)$ and $(4^4,2)$, then $\pi_2$ is graphic by lemma
2.3 and lemma 2.4. Let $G_1$ be a realization of $\pi_2$, then
$(K_6-2K_2)\cup G_1$ is a realization of
$\pi=(5^2,4^k,3^t,2^{n-2-k-t})$. Since $K_6-C_4\subseteq K_6-2K_2$,
$\pi$ is potentially $K_6-C_4$-graphic. If $n=11$ and
$\pi_2=(4^3,2^2)$ or $(4^4,2)$, then $\pi=(5^2,4^7,2^2)$ or
$(5^2,4^8,2)$. If $n\leq10$, then $\pi$ is one of the following:
$(5^2,4^4,2)$, $(5^2,4^4,2^2)$, $(5^2,4^4,2^3)$, $(5^2,4^4,2^4)$,
$(5^2,4^4,3^2,2)$,  $(5^2,4^4,3^2,2^2)$, $(5^2,4^5,2)$,
$(5^2,4^5,2^2)$, $(5^2,4^5,2^3)$, $(5^2,4^5,3^2,2)$, $(5^2,4^6,2)$,
$(5^2,4^6,2^2)$, $(5^2,4^7,2)$. It is easy to check that all of
these are potentially $K_6-C_4$-graphic.
\par
\textbf{Subcase 2.2:} $i=3$. Then $\pi=(5^3,4^k,3^t,2^{n-3-k-t})$.
If $k=0$, then $\pi=(5^3,3^t,2^{n-3-t})$. Since
$\pi\neq(5^3,3^3,2)$, we have $n\geq8$. If $n\geq9$, it is enough to
show $\pi_1=(3^{t-3},2^{n-2-t})$ is graphic. It follows by lemma
2.3. If $n=8$, then $\pi=(5^3,3^3,2^2)$ which is potentially
$K_6-C_4$-graphic.
\par
If $k=1$, then $\pi=(5^3,4,3^t,2^{n-4-t})$. Let
$\pi_1=(5^3,4,3^3,2)$, $\pi_2=(3^{t-3},2^{n-5-t})$. It is easy to
see that $\pi_1$ is potentially $K_6-C_4$-graphic. Let $G_1$ be a
realization of $\pi_1$ with $K_6-C_4\subseteq G_1$. If $n\geq12$,
then $\pi_2$ is graphic by lemma 2.3. Let $G_2$ be a realization of
$\pi_2$, then $G_1\cup G_2$ is a realization of $\pi$. If $n\leq11$,
then $\pi$ is one of the following: $(5^3,4,3^3,2)$,
$(5^3,4,3^3,2^2)$, $(5^3,4,3^3,2^3)$, $(5^3,4,3^3,2^4)$,
$(5^3,4,3^5,2)$, $(5^3,4,3^5,2^2)$. It is easy to check that all of
these are potentially $K_6-C_4$-graphic. Similarly, one can show
that $\pi=(5^3,4^k,3^t,2^{n-3-k-t})$ is potentially
$K_6-C_4$-graphic for the cases $k=2$ and $k=3$.
\par
Now we consider the case where $k\geq4$. If $n\geq12$, let
$\pi_1=(5^3,4^3,3)$, $\pi_2=(4^{k-3},3^{t-1},2^{n-3-k-t})$. It is
easy to see that $\pi_1$ is potentially $K_6-C_4$-graphic. Let $G_1$
be a realization of $\pi_1$ with $K_6-C_4\subseteq G_1$. If
$\pi_2\neq(4^3,2^2)$ and $(4^4,2)$, then $\pi_2$ is graphic by lemma
2.3 and lemma 2.4. Let $G_2$ be a realization of $\pi_2$, then
$G_1\cup G_2$ is a realization of $\pi=(5^3,4^k,3^t,2^{n-3-k-t})$.
If $\pi_2=(4^3,2^2)$ or $(4^4,2)$, then $\pi=(5^3,4^6,3,2^2)$ or
$(5^3,4^7,3,2)$. If $n\leq11$, then $\pi$ is one of the following:
$(5^3,4^4,3,2)$, $(5^3,4^4,3,2^2)$, $(5^3,4^4,3,2^3)$,
$(5^3,4^4,3^3,2)$, $(5^3,4^5,3,2)$, $(5^3,4^5,3,2^2)$,
$(5^3,4^6,3,2)$. It is easy to check that all of these are
potentially $K_6-C_4$-graphic.
\par
If $\pi^\prime$ does not satisfy $(2)$, then
$\pi^\prime=(d_1^\prime,d_2^\prime,d_3^\prime,3^k,2^{n-4-k})$ and
$d_1^\prime+d_2^\prime+d_3^\prime> n-1+2k+n-4-k+1$, i.e.,
$d_1^\prime+d_2^\prime+d_3^\prime>2n+k-4 $. Hence,
$\pi=(d_1,d_2,d_3,3^k,2^{n-3-k})$ and $d_1+d_2+d_3> 2n+k-2$, a
contradiction.
\par
If $\pi^\prime$ does not satisfy $(3)$, then $\pi^\prime$ is one of
the following: $(5^2,4^6)$, $(5^2,4^7)$, $(6^2,3^6)$, $(6,5,4,3^5)$,
  $(6,5,3^7)$, $(5^3,4,3^3)$, $(5^3,3^5)$, $(5^2,4^2,3^4)$,
  $(5^2,4,3^6)$, $(5^2,4,3^4)$, $(5^2,3^6)$, $(6,5,3^5,2)$,
  $(5^3,3^3,2)$, $(5^2,4,3^4,2)$, $(5^2,3^6,2)$, $(5^2,3^4,2)$,
  $(5^2,3^4,2^2)$, $(6,5,3^5)$, $(7,5,3^6)$. Hence, $\pi$ is one of
the following: $(6^2,4^6,2)$, $(6,5^2,4^5,2)$, \ \ $(5^4,4^4,2)$, \
\ $(6^2,4^7,2)$, \ \ $(6,5^2,4^6,2)$, \ \ $(5^4,4^5,2)$, \ \
$(7^2,3^6,2)$, \ \ $(7,6,4,3^5,2)$,$(7,5^2,3^5,2)$,$(7,6,3^7,2)$,
$(6^2,5,4,3^3,2)$,$(6,5^3,3^3,2)$,$(6^2,5,3^5,2)$,
$(6^2,4^2,3^4,2)$,$(6,5^2,4,3^4,2)$,$(5^4,3^4,2)$,
$(6^2,4,3^6,2)$,$(6,5^2,3^6,2)$,$(6^2,4,3^4,2)$, $(6,5^2,3^4,2)$,
$(6^2,3^6,2)$, $(7,6,3^5,2^2)$, $(6^2,5,3^3,2^2)$,
$(6^2,4,3^4,2^2)$, $(6,5^2,3^4,2^2)$, $(6^2,3^6,2^2)$,
$(6^2,3^4,2^2)$, $(6^2,3^4,2^3)$, $(7,6,3^5,2)$, $(8,6,3^6,2)$. It
is easy to check that all of these are potentially
$K_6-C_4$-graphic.
\par
\textbf{Case 4:} $d_n=1$. Consider
$\pi^\prime=(d_1^\prime,d_2^\prime,\cdots,d_{n-1}^\prime)$ where
$d_1^\prime\geq5$, $d_2^\prime\geq4$  and $d_6^\prime\geq3$. If
$\pi^\prime$ satisfies $(1)$-$(3)$, then by the induction
hypothesis, $\pi^\prime$ is potentially $K_6-C_4$-graphic, and hence
so is $\pi$.
\par
  If $\pi^\prime$ does not satisfy $(1)$, i.e., $d_2^\prime=4$, then $\pi=(5^2,4^k,3^t,2^i,1^{n-2-k-t-i})$ where
$k+t\geq4$, $n-2-k-t-i\geq1$ and $n-k-i$ is even. We will show that
$\pi$ is potentially $K_6-C_4$-graphic.
\par
If $k=0$, then $\pi=(5^2,3^t,2^i,1^{n-2-t-i})$. If $n\geq10$, we
only need to show that $\pi_1=(3^{t-4},2^i,1^{n-2-t-i})$ is graphic.
Since $\pi\neq(5^2,3^7,1)$ and $(5^2,3^6,1^2)$, then
$\pi_1\neq(3^3,1)$, $(3^2,1^2)$.  By lemma 2.3, $\pi_1$ is graphic.
 If $n\leq9$, since $\pi\neq(5^2,3^5,1)$ and $(5^2,3^5,2,1)$,
then $\pi=(5^2,3^4,1^2)$ or $(5^2,3^4,2,1^2)$. It is easy to check
that both of them are potentially $K_6-C_4$-graphic.
\par
 If $k=1$,
then $\pi=(5^2,4,3^t,2^i,1^{n-3-t-i})$. If $n\geq9$, it is enough to
show $\pi_1=(3^{t-3},2^i,1^{n-2-t-i})$ is graphic. Since
$\pi\neq(5^2,4,3^6)$ and $(5^2,4,3^5,1)$, then $\pi_1\neq(3^3,1)$,
$(3^2,1^2)$.  By lemma 2.3, $\pi_1$ is graphic.
 If $n\leq8$,
then $\pi=(5^2,4,3^3,1)$ or $(5^2,4,3^3,2,1)$. It is easy to check
that both of them are potentially $K_6-C_4$-graphic. With the same
argument as above, one can show that
$\pi=(5^2,4^k,3^t,2^i,1^{n-2-k-t-i})$ is potentially
$K_6-C_4$-graphic for the cases $k=2$ and $k=3$.
\par
Now we consider the case where $k\geq4$. If $n\geq11$, let
$\pi_1=(5^2,4^4)$, $\pi_2=(4^{k-4},3^{t},2^i,1^{n-2-k-t-i})$.  If
$\pi_2\neq(4,3^2,1^2)$, $(4,3,1^3)$, $(4^2,2,1^2)$, $(4^2,3,2,1)$,
$(4^3,1^2)$, $(4^3,3,1)$, $(4^2,3,1^3)$, $(4^2,1^4)$, $(4^3,2,1^2)$,
$(4^4,1^2)$ and $(4^3,1^4)$, then $\pi_2$ is graphic by lemma 2.3
and lemma 2.4. Let $G_1$ be a realization of $\pi_2$, then
$(K_6-2K_2)\cup G_1$ is a realization of
$\pi=(5^2,4^k,3^t,2^i,1^{n-2-k-t-i})$. Since $K_6-C_4\subseteq
K_6-2K_2$, $\pi$ is potentially $K_6-C_4$-graphic. If $\pi_2$ is one
of the following: $(4,3^2,1^2)$, $(4,3,1^3)$, $(4^2,2,1^2)$,
$(4^2,3,2,1)$, $(4^3,1^2)$, $(4^3,3,1)$, $(4^2,3,1^3)$, $(4^2,1^4)$,
$(4^3,2,1^2)$, $(4^4,1^2)$, $(4^3,1^4)$, then $\pi$  is one of the
following: $(5^2,4^5,3^2,1^2)$, $(5^2,4^5,3,1^3)$,
$(5^2,4^6,2,1^2)$, $(5^2,4^6,3,2,1)$, $(5^2,4^7,1^2)$,
$(5^2,4^7,3,1)$, $(5^2,4^6,3,1^3)$, $(5^2,4^6,1^4)$,
$(5^2,4^7,2,1^2)$, $(5^2,4^8,1^2)$, $(5^2,4^7,1^4)$. If \ $n\leq10$,
\ then \ $\pi$ \ is \ one \ of \ the \ following: $(5^2,4^4,1^2)$,
$(5^2,4^4,2,1^2)$, $(5^2,4^4,2^2,1^2)$, $(5^2,4^4,3^2,1^2)$, \ \
$(5^2,4^4,3,1)$, \ \  $(5^2,4^4,3,2,1)$, \ \  $(5^2,4^4,3,2^2,1)$, \
\  $(5^2,4^4,3^3,1)$, \ \ $(5^2,4^4,3,1^3)$, \ \ $(5^2,4^4,1^4)$,
$(5^2,4^5,1^2)$, \  $(5^2,4^5,2,1^2)$, $(5^2,4^5,3,1)$,
$(5^2,4^5,3,2,1)$, $(5^2,4^6,1^2)$, $(5^2,4^6,3,1)$. It is easy to
check that all of these are potentially $K_6-C_4$-graphic.
\par
If $\pi^\prime$ does not satisfy $(2)$, then
$\pi^\prime=(d_1^\prime,d_2^\prime,d_3^\prime,3^k,2^t,1^{n-4-k-t})$
and $d_1^\prime+d_2^\prime+d_3^\prime> n-1+2k+t+1$, i.e.,
$d_1^\prime+d_2^\prime+d_3^\prime>n+2k+t $. Hence,
$\pi=(d_1,d_2,d_3,3^k,2^t,1^{n-3-k-t})$ and $d_1+d_2+d_3> n+2k+t+1$,
a contradiction.
\par
If $\pi^\prime$ does not satisfy $(3)$, since $\pi\neq(6,5,3^5,1)$
and $(n-1,5,3^6,1^{n-8})$, then $\pi^\prime$ is one of the
following: $(5^2,4^6)$, $(5^2,4^7)$, $(6^2,3^6)$, $(6,5,4,3^5)$,
  $(6,5,3^7)$, $(5^3,4,3^3)$, $(5^3,3^5)$, $(5^2,4^2,3^4)$,
  $(5^2,4,3^6)$, $(5^2,4,3^4)$,  $(6,5,3^5,2)$,
  $(5^3,3^3,2)$, $(5^2,4,3^4,2)$, $(5^2,3^6,2)$,  \ $(5^2,3^4,2)$,
  $(5^2,3^4,2^2)$, \ $(6,5,3^6,1)$, \ $(5^3,3^4,1)$, \ $(5^2,4,3^5,1)$,
  $(5^2,3^7,1)$, $(5^2,3^6,1^2)$, $(5^2,3^5,1)$,
   $(6,5,3^5)$. Since $\pi\neq(5^3,3^4,1)$ and $(n-1,5,3^5,1^{n-7})$,
   then $\pi$ is one of the following:
   $(6,5,4^6,1)$, $(5^3,4^5,1)$, $(6,5,4^7,1)$, $(5^3,4^6,1)$, \ \
   $(7,6,3^6,1)$,\ \
    $(7,5,4,3^5,1)$,\ \  $(6^2,4,3^5,1)$,\ \  $(7,5,3^7,1)$,\ \  $(6^2,3^7,1)$,\ \
    $(6,5^2,4,3^3,1)$,
    $(5^4,3^3,1)$,$(6,5^2,3^5,1)$,$(6,5,4^2,3^4,1)$,
    $(5^3,4,3^4,1)$,  $(6,5,4,3^6,1)$,$(5^3,3^6,1)$,
    $(6,5,4,3^4,1)$,$(7,5,3^5,2,1)$,$(6^2,3^5,2,1)$,
    $(6,5^2,3^3,2,1)$, $(6,5,4,3^4,2,1)$, \  $(5^3,3^4,2,1)$, \
    $(6,5,3^6,2,1)$, \  $(6,5,3^4,2,1)$, \  $(6,5,3^4,2^2,1)$, \
    $(7,5,3^6,1^2)$, $(6^2,3^6,1^2)$, $(6,5^2,3^4,1^2)$,
    $(6,5,4,3^5,1^2)$, $(5^3,3^5,1^2)$, $(6,5,3^7,1^2)$,
    $(6,5,3^6,1^3)$, $(6,5,3^5,1^2)$, $(6^2,3^5,1)$. It is easy to check that all of
these are potentially $K_6-C_4$-graphic.

\par
\vspace{0.5cm}

\section{  Application }

\par
In the remaining of this section, we will use theorem 3.1 to find
exact values of $\sigma(K_6-C_4,n)$. Note that the value of
$\sigma(K_6-C_4,n)$ was determined by Hu and Lai in $[7]$ so a much
simpler proof is given here.
\par
  \textbf{Theorem }  (Hu and Lai [7])  If $n\geq6$, then
    $\sigma(K_6-C_4,n)=6n-10$.
\par
\textbf{Proof:} First we claim that for $n\geq6$,
$\sigma(K_6-C_4,n)\geq 6n-10$. Take $\pi_1=((n-1)^3,3^{n-3})$, then
$\sigma(\pi_1)=6n-12$, and it is easy to see that $\pi_1$ is not
potentially $K_6-C_4$-graphic by condition (2) in Theorem 3.1.
\par
  Now we show that if $\pi$ is an $n$-term $(n\geq6)$ graphic sequence
with $\sigma(\pi) \geq 6n-10$, then there exists a realization of
$\pi$ containing a $K_6-C_4$.
\par
If $d_2\leq4$, then $\sigma(\pi)\leq d_1+4(n-1)\leq
n-1+4(n-1)=5n-5<6n-10$, a contradiction. Hence, $d_2\geq5$.
\par
 If $d_6\leq2$, then $\sigma(\pi)\leq
d_1+d_2+d_3+d_4+d_5+2(n-5)\leq20+2(n-5)+2(n-5)=4n<6n-10$, a
contradiction. Hence, $d_6\geq3$.
\par
  Since \ \  $\sigma(\pi)\geq6n-10$, \ \ then \ \ $\pi$ \ \ is \ \   not \ \  one \ \   of
the \ \ following:
 $(d_1,d_2,d_3,3^k,2^t,1^{n-3-k-t})$, $(5^2,4^6)$, $(5^2,4^7)$, $(6^2,3^6)$, $(6,5,4,3^5)$,
  $(6,5,3^7)$, $(5^3,4,3^3)$, $(5^3,3^5)$, $(5^2,4^2,3^4)$,
  $(5^2,4,3^6)$, $(5^2,4,3^4)$, $(5^2,3^6)$, $(6,5,3^5,2)$,
  $(5^3,3^3,2)$, $(5^2,4,3^4,2)$, $(5^2,3^6,2)$,$(5^2,3^4,2)$,
  $(5^2,3^4,2^2)$,$(6,5,3^6,1)$,$(5^3,3^4,1)$, $(5^2,4,3^5,1)$,
  $(5^2,3^7,1)$,  $(5^2,3^6,1^2)$,  $(5^2,3^5,1)$, \ \
   $(n-1,5,3^5,1^{n-7})$, \ \  $(n-1,5,3^6,1^{n-8})$.
 Thus, $\pi$ satisfies the conditions (1)-(3) in
Theorem 3.1. Therefore, $\pi$ is potentially $K_6-C_4$-graphic.

\par

\end{document}